\documentclass{article}

\usepackage{amsmath,amssymb,amsfonts,bbm,color,bm,nicefrac,
graphicx,enumerate}
\usepackage{amsthm}
\usepackage{epstopdf}
\usepackage{geometry}
\usepackage{caption}
\usepackage[toc,page,title]{appendix}
\usepackage{authblk}

\newtheorem{thm}{Theorem}[section]
\newtheorem{lem}[thm]{Lemma}
\newtheorem{prop}[thm]{Proposition}
\newtheorem{cor}[thm]{Corollary}

\newtheorem{remark}[thm]{Remark}

\newcommand{\comment}[1]{}

\newcommand{\tmtextbf}[1]{{\bfseries{#1}}}

\begin{document}
\nocite{*}
\title{Coarsening with a Frozen Vertex}

\author{Michael Damron$^{1,}$\thanks{Partially supported by US-NSF Grant DMS-1419230},
       %Sinziana M. Datcu$^{2,}$\thanks{Partially supported by NSF Grants, 0ISE-0730136 and DMS-1007524.},
       Hana Kogan$^{2}$, \\ %\thanks{}, \\
       Charles M. Newman$^{3,4,}$\thanks{Partially supported by US-NSF Grants 
                                         DMS-1007524 and DMS-1507019},
       Vladas Sidoravicius$^{3,4}$ %\thanks{ADD GRANT ACKNOWLEDGMENT?}
       \newline
       \\
       $^1$ Georgia Institute of Technology, Atlanta GA and 
       Indiana University, Bloomington IN \\
       $^2$ CUNY College of Staten Island \\
       $^3$ Courant Institute of Mathematical Sciences \\
       $^4$ NYU--ECNU Institute of Mathematical Sciences at NYU Shanghai \\
       %$^5$ Department of Mathematics, University of California, Irvine \\
       %$^6$ IMPA
       }
\maketitle

\begin{abstract}
\noindent
In the standard nearest-neighbor coarsening model with state space
$\{-1,+1\}^{\mathbb{Z}^2}$ and initial state chosen from symmetric product
measure, it is known (see~\cite{NNS}) that almost surely, every vertex
flips infinitely often. In this paper, we study the modified model in which
a single vertex is frozen to $+1$ for all time, and show that 
every other site still flips infinitely often. The proof combines stochastic
domination (attractivity) and influence propagation arguments.

\end{abstract}

\section{Introduction}

As in our earlier paper \cite{DDKNS}, we study and compare the long time behavior of 
two continuous time Markov coarsening models with state space $\Omega=\{-1,+1\}^{\mathbb{Z}^d}$. 
One, $\sigma(t)$, is the standard model in which at time zero 
$\{\sigma_x(0):x\in \mathbb{Z}^d\}$ is an i.i.d. set with $\theta \equiv P(\sigma_x(0)=+1)=1/2$ 
and then vertices update to agree with a strict majority of their $2d$ nearest 
neighbors or, in case of a tie, choose their value by tossing a fair coin. The 
modified model, $\sigma'(t)$, is the same except that $\sigma'$ at the origin 
$(0,0....0)$ is frozen to $+1$ for all $t\geq 0$.

For $d=2$, it is an old result \cite{NNS} that in the standard $\sigma(t)$  model, almost surely, every 
vertex changes sign infinitely many times as $t \to \infty$. The main 
result of this paper (see Theorem  \ref{5}) is that the same is true for the frozen 
model $\sigma'(t)$ on $\mathbb{Z}^2$. It is believed 
(see, for example, Sec.~6.2 of~\cite{NS}), % \cite{FSS}, 
but not proved, that the 
$d=2$ behavior of $\sigma$ remains valid at least for some values of $d>2$. If this were so,
then the arguments of this paper would show the same for the corresponding $\sigma'$ model. 

In the previous paper \cite{DDKNS} we  considered models with infinitely many frozen vertices
and in this paper a model with a single frozen vertex. 
It would be of interest to study models with finitely many, but more than one, frozen vertices; 
in this regard, see the remark following the proof of Theorem~\ref{four} below.

\section{Results}
In this section we fix $d=2$. We also use the standard convention that the updates are made when independent rate one Poisson process clocks at each vertex ring. 

Let $A_T$ denote the event that the ``right" neighbor of the origin 
(at $x=(1,0)$) is $-1$ 
%(i.e., $\sigma'_{(1,0)}(t) =-1$ )  
for some $t \geq T$.
Let $A_T^1 \subset A_T$ denote the event that the right neighbor of the 
origin is the { \it first} neighbor to be $-1$ at some time $t \geq T$
(more precisely, that no other neighbor is $-1$ at an earlier time in $[T,\infty]$). 
Let $B_{L,s}$ for $s \in \{-1, +1\}^{\Lambda _L}$ (where $\Lambda_L=\{-L, -L+1, ....,L\}^2$) 
denote the event that $\sigma'(0)|_{\Lambda_L}=s$ and write $B^{ }_{L,+}$ when $s \equiv +1$.
We denote the probability measure for the frozen origin $\sigma'(\cdot)$ model by $P'$ and 
that for the regular coarsening model $\sigma(\cdot)$ by~$P$.

\begin{lem} \label{L1} For all $L$,
$$P(A_0^1| B_{L, +})\geq 1/4 \, .$$

\end{lem}

\proof 

The result is an easy consequence of symmetry among the four neighbors of  the
origin and the fact that $P(A_0) = 1$
(indeed, for all $T$, $P(A_T) = 1$ --- see \cite{NNS}). 
%$$P(A_0  \hspace{.1cm} \tt{i.o.})=1$ a.s.
\qed

Let $\Sigma_T^L$ denote the sigma-field generated by  the initial spin values
and clock rings and coin tosses up to time $T$ inside the box $\Lambda_L$.
\begin{prop}\label{P1}
For any $T$, $L$, 

$$P'(A_T|\Sigma^L_T) \geq 1/4\hspace{.1cm}\mathrm{a.s.} \, . $$
%for all $T$ and $L$.
\end{prop}

\proof

Let $\tilde \sigma_T^L(.)$ denote the model with the spin values at all sites in $\Lambda_L$ frozen 
to $+1$ from time $0$ up to time $T$ and with the spin value  at the origin remaining frozen 
at $+1$ thereafter. Denote the corresponding probability
measure by $\tilde P_T^L$.
Under the standard coupling,  $\tilde \sigma(\cdot)$ stochastically dominates $\sigma'(\cdot)$, 
so we have $$P'(A_T|\Sigma_T^L)\geq \tilde P_T^L(A_T)\geq \tilde P_T^L(A_T^1).$$

To continue the proof, we will use the following result about the ``propagation speed" of
influence between different spatial regions: 

\begin{lem}\label{finprop}
Let $D_T^L$ denote the event that  $\sigma_x(t) = +1 \, \forall x \in \Lambda_L, \forall  t \in [0,T]$.
Then 
$$ \forall L,T,\varepsilon, \exists L' \ \mathrm{such } \  \mathrm{that}\   P(D_T^L | B_{L',+}) \geq 1-\varepsilon \, .$$

\end{lem}

\proof  %of Lemma \ref{finprop}. 

Let $L'>>L$ and note that  given $B_{L',+}$, $(D_T^L)^c$ can occur only if there is a nearest 
neighbor (self-avoiding) path between the boundaries 
of the two sets, $\mathbb{Z}^2\setminus \Lambda_{L'}$ and $\Lambda_{L}$,
along which there are clock rings occurring in succession
between times $0$ and $T$. Any such path is at least of length $L'-L$
(i.e., contains at least $L'-L$ vertices besides the starting one).

Consider a particular path $\gamma$  of length $m\geq L'-L$. 
For each $m$ there are no more than $3^m$ such paths from each 
boundary point and the  time it takes for successive clock rings along $\gamma$ is 
at least $S_m=\sum_{i=1}^m \tau_i$ where the $\tau_i$ are i.i.d.
exponential random variables with parameter 1. 
By the exponential Markov inequality,  for any $\alpha >0$,
$$P(\sum_{i=1}^m \tau_i<T)=P(-\sum_{i=1}^m \tau_i>-T) \leq 
\frac{E(e^{-\alpha\sum_{i=1}^m \tau_i})}{E(e^{-\alpha T})} = 
e^{\alpha T}E\{e^{-\alpha \tau_i}\}^m=\frac{e^{\alpha T}}{(1+\alpha)^m} .$$

Therefore, since there are at most $CL'$ possible starting points (for some constant $C$), 
$$P( (D_T^L)^c|B_{L',+})\leq C\,L'\sum_{m=L'-L}^\infty 3^m 
\frac{e^{\alpha T}}{(1+\alpha)^m}= C(\alpha, T, L)L' (\frac{3}{1+ \alpha})^{L'} \, ,$$
where $C(\alpha, T, L)$ is a constant depending on $\alpha$, $T$ and $L$. Taking 
$\alpha >2$ and the limit as $L' \to \infty$ completes the proof of the lemma.

\proof We continue the proof of Proposition~\ref{P1}. 
Pick $\epsilon>0$ and fix $T$ and $L$.  By Lemma \ref{finprop},  $\exists$ $L'$ such that 
$$ P( D_T^L|B_{L',+}) \geq 1-\epsilon \, .$$
% $$P(\sigma_t(x)=+1 \, \forall x \in B_L,\forall  t \in [0,T]\, |\, B_{L',+}) \geq 1-\epsilon \, .$$ 
%  for all $x \in B_L$. 
Therefore, given $B_{L',+}$,   with probability at least $1-\epsilon$, 
$\sigma_t(\cdot)$ positively dominates $\tilde \sigma_L^T(\cdot)$   
for $0\leq t < S$, where $S =\inf \{t>0\, | \, \sigma_t(0,0)=-1\}$,  %\{t>0,A_T|B^0_{L',+}\}$ in $\sigma()$ )  
and so 

$$\tilde{P}_L^T(A_T^1) \geq P(A_T^1|B_{L',+})-\epsilon\geq 1/4-\epsilon \, .$$
Taking the limit as $\epsilon \to 0$ completes the proof of Proposition~\ref{P1}.
%$ $\lim_{\epsilon \to 0}P_L^T(A_T^1) \geq 1/4$ obtain the Proposition \ref {one}.

\qed

Now let $\Sigma_T$ denote the sigma field generated by the initial assignment 
of spins on $\mathbb{Z}^2$ and the clock rings 
and coin tosses on $\mathbb{Z}^2$ up to time $T$. 

\begin{prop} \label{two}
For all $T$,
$$P'(A_T|\Sigma_T) \geq 1/4\hspace{.1cm} \mathrm{a.s.} \, .$$ 
\end{prop}

\proof
For $L\geq 1$ let $X_L= P'(A_T|\Sigma^L_T)$. $\{ \Sigma^L_T, L\geq 1\}$ is 
an increasing filtration of sigma fields, and  
$E(X_{L+1}|\Sigma^L_T )=X_L$.
%$E(X_{L+1}|X_L)=X_L$.
By the martingale convergence theorem, $\lim_{L\to \infty} (X_L)=X_\infty=P'(A_T|\Sigma_T)$ 
and since $X_L\geq 1/4$ for all $L$, we have $P'(A_T|\Sigma_T) \geq 1/4$.

\qed

Let $A_{T,T'}$ denote the event that the right neighbor of the origin is $-1$ for
some time $t \in [T,T']$. The following is immediate from Proposition~\ref{two}.

\begin{cor} \label{three}
$$\lim_{T'\to \infty}P'(A_{T,T'}|\Sigma_T) \geq 1/4\hspace{.1cm}\mathrm{a.s.}\, .$$
\end{cor}
%follows immediately from Proposition~\ref{two}.

\begin{lem}\label{4}
For any $T \geq 0$ and $ \gamma >0$, $\exists$ a deterministic $T'$ such that  
$$P\{ \omega:P'(A_{T,T'}|\Sigma_T)\geq 1/8\}\geq 1-\gamma \, .$$
\end{lem}

\proof
This is a straightforward consequence of the preceding corollary.
%For a.e. $\omega \in \Omega$, $\exists$ $ T_\omega<\infty $ so that 
%$A_{T, T_\omega}$ occurs.  $T_\omega$ is a random time and we  let 
%$T_M=\min \{ T_\omega, M\}$. Then $P'(A_{T, T_M}|\Sigma_T) $ is increasing in $M$. 
%Choosing $T'$ large enough yields the lemma. 
\qed

\begin{thm}\label{5}

For any $T$, 
$$ P'(A_T) = 1, \mathrm{and} \, \mathrm{ hence }\, P'(\cap_{T>0} A_T) =1 \, .$$
%$$P\{ \omega: P'(A _T\hspace{.1cm} \mathrm{i.o.})=1\}=1$$
It follows that with probability one, $\sigma'_{(1,0)}(t)$ changes sign infinitely
many times as $t \to \infty$.
\end{thm}

\proof

Given $T$ and $\epsilon>0$ construct a sequence of deterministic times $\{T_i; i\geq 0\}$  so that

1. $T_0=T$, and

2. $P'\{\omega: P'(A_{T_{i-1}, T_i}|\Sigma_{T_{i-1}}) \geq 1/8\}\geq 1-\frac{\epsilon}{2^i}$.

%Therefore, except for an event of probability at most  
%$ \sum_{i=1}^\infty\frac{\epsilon}{2^i}=\epsilon$, the sequence  of 
%occurrences of $A_{T_i, T_{i+1}}$ dominates the successes of i.i.d. trials  
%with success probability $1/8$ and so 
Condition now on the event (of probability at least 
$1 - \sum_{i=1}^\infty\frac{\epsilon}{2^i}= 1- \epsilon$) that
$P'(A_{T_{i-1}, T_i}|\Sigma_{T_{i-1}}) \geq 1/8$ for all $i$. On this conditioned
probability space, letting $\tilde{W}_i =1$ (and otherwise~$0$) if $A_{T_{i-1}, T_i}$ occurs,
we note that the $\tilde{W}_i$'s stochastically dominate i.i.d. $\{0,1\}$-valued $W_i$'s
with Prob$(Wi=1)= 1/8$. Thus
 
%$$P\{\omega;P'( A_{T_i, T_{i+1}}\hspace{.1cm} 
%\mathrm {finitely \hspace{.1cm}often})=\lim_{n \to \infty} (7/8)^n=0\}
%$$P'(A_{T_i, T_{i+1}} \mathrm{occurs }\, \mathrm{for }\, \mathrm{only }\, \mathrm{finitely }\,
%\mathrm{many }\,  i)
%\leq \sum_{i=1}^\infty\frac{\epsilon}{2^i}=\epsilon \, .$$
$$P'(A_{T_{i-1}, T_i} \mathrm{occurs }\, \mathrm{for }\, \mathrm{only }\, \mathrm{finitely }\,
\mathrm{many }\,  i) \leq \epsilon \, .$$
Letting $\epsilon \to  0$ completes the proof of the first part of the theorem. The
second part then follows because by stochastic domination (attractivity) and the results
of~\cite{NNS}, $\sigma'_{(0,0)}(t_i)$ equals $+1$ for an infinite sequence of $t_i \to 
\infty$. 
\qed

The next theorem follows from 
a modified version of the proof of Theorem~\ref{5}. 
\begin{thm}\label{four}
Every site in $\{\mathbb{Z}^2\setminus (0,0)\}$  flips infinitely many times in $\sigma'(\cdot)$
with probability one. 
\end{thm}
\proof 

For any site $z$ other than the origin,  and for $L$ much larger than say the Euclidean norm
of $z$, we consider the unfrozen $\sigma$ model in which at time zero all the vertex values
are set to $+1$
in the box of side length $2L$, centered at $z/2$ (so that the origin and z are located
symmetrically with respect to this box). Then with probability $1/2$ the vertex at $z$
flips to $-1$ before the one at the origin flips and until just after that time, there is no
difference between the frozen (at the origin) $\sigma'$ model and the unfrozen $\sigma$ model.
Hence there is probability at least $1/2$ in $\sigma'$ that $z$ will flip to minus. 
By applying the methods used in the proof of Theorem~\ref{5} (but with $1/4$ now
replaced by $1/2$), we conclude that $z$ will flip infinitely many times with
probability one. 
\qed

We note that the line of reasoning in the proof of the last theorem
could have also been used to give a modified proof of 
Theorem~\ref{5} with $1/4$ replaced by $1/2$.
A more interesting remark is the following.
\begin{remark}\label{Remark1}
For the process $\sigma''$ with some finite set ${\cal S}$ of vertices frozen to
$+1$, it is possible to show by an extension of the arguments used in this paper that
there is a finite deterministic ${\cal S}' \supseteq {\cal S}$ such that all sites
in $\{\mathbb{Z}^2\setminus {\cal S}' \}$ flip infinitely many times in $\sigma''(\cdot)$
with probability one. In some cases,  ${\cal S}'$ must be strictly larger than ${\cal S}$ --- e.g.,
when ${\cal S} =  \{(-L,-L), (-L+L), (+L,-L), (+L,+L) \}$, $ {\cal S}'$ includes all of $\Lambda_L$.  
One may also consider processes where some vertices  are frozen to $-1$ 
and some to $+1$. 
We expect to  to pursue these issues in a future paper.

\end{remark}
%\section{Frozen box.} 

\end{document}